\documentclass[11pt]{article}
\usepackage[margin=1.1in]{geometry}
\usepackage{amsmath,amssymb,amsthm}
\usepackage{enumitem}
\usepackage{tikz}
\usepackage{xcolor}
\usepackage[colorlinks=true,linkcolor=blue,citecolor=blue,urlcolor=blue]{hyperref}

\newtheorem{theorem}{Theorem}[section]
\newtheorem{lemma}[theorem]{Lemma}
\newtheorem{proposition}[theorem]{Proposition}
\newtheorem{corollary}[theorem]{Corollary}
\newtheorem{conjecture}[theorem]{Conjecture}
\theoremstyle{definition}
\newtheorem{question}[theorem]{Question}

\newcommand{\R}{\mathbb{R}}
\newcommand{\Z}{\mathbb{Z}}
\newcommand{\Q}{\mathbb{Q}}

\begin{document}

\title{The Gerards--Schrijver conjecture\\ on the Edmonds--Johnson property is false}

\author{Alberto Del Pia\\[2pt]
{\small Department of Industrial and Systems Engineering \& Wisconsin Institute for Discovery}\\
{\small University of Wisconsin--Madison, Madison, WI, USA}\\
{\small\texttt{delpia@wisc.edu}}}

\date{November 20, 2026}

\maketitle

\begin{abstract}
An integral matrix $A$ has the Edmonds--Johnson property if, for all integral vectors $b,c,l,u$, the
first Chv\'atal closure of $\{x: b\le Ax\le c,\ l\le x\le u\}$ equals its integer hull. Gerards and
Schrijver conjectured that a totally half-modular matrix has the Edmonds--Johnson property if and only
if it has no minor equal to one of two specific matrices, of size $3\times3$ and $3\times4$. We disprove
this conjecture. We give a totally half-modular $5\times6$ matrix that does not have the Edmonds--Johnson
property and has no minor equal to either matrix. It is a minimal obstruction: every proper minor has the
Edmonds--Johnson property. The proofs work on a graph whose edges are the rows and bounds of the matrix, and
on which the minor operations act by deleting, contracting and halving edges. For minimality, we reduce every proper minor to
two classes of matrices known to have the Edmonds--Johnson property.
\end{abstract}

\section{Introduction}\label{sec:intro}

For a rational polyhedron $P\subseteq\R^n$, the
\emph{integer hull} of $P$ is $P_I:=\operatorname{conv}(P\cap\Z^n)$. A \emph{Chv\'atal--Gomory cut}
(\emph{CG cut})~\cite{Gom58,Chv73} for $P$ is an inequality $wx\le\lfloor\delta\rfloor$, where $w\in\Z^n$ and $wx\le\delta$
is valid for $P$. The \emph{Chv\'atal closure} $P'$ of $P$ is the set of points of $P$ that satisfy
all CG cuts for $P$. Chv\'atal~\cite{Chv73} and Schrijver~\cite{Sch80} showed that $P'$ is a rational
polyhedron and that $P_I$ is obtained from $P$ by finitely many applications of this operation. The
\emph{Chv\'atal rank} of $P$ is the minimum number of applications needed.

Unlike the Chv\'atal rank, which depends on the polyhedron, the following notion of rank depends only on the constraint matrix. The \emph{strong Chv\'atal rank} of an
integral matrix $A$ is the maximum Chv\'atal rank of the polyhedra
\begin{equation}\label{eq:system}
P(A;b,c,l,u):=\{x\in\R^n:\ b\le Ax\le c,\ l\le x\le u\}
\end{equation}
over all integral vectors $b,c,l,u$. By a theorem of Cook, Gerards, Schrijver and
Tardos~\cite{CooGerSchTar86}, this maximum exists. The matrices of strong Chv\'atal rank $0$ are exactly the totally unimodular
matrices (see~\cite[Chapter~19]{Sch86}). Following~\cite{GerSch86,dPMusZam18JCTB}, $A$ has the \emph{Edmonds--Johnson property}
(\emph{EJ property}) if its strong Chv\'atal rank is at most $1$, that is, if $P'=P_I$ for every polyhedron
$P$ of the form~\eqref{eq:system}.

Few classes of matrices with the EJ property are known. Edmonds and Johnson~\cite{EdmJoh73} showed that
every integral matrix whose columns have $\ell_1$-norm at most $2$ has the EJ property. Gerards and
Schrijver~\cite{GerSch86} proved that an integral matrix whose rows have $\ell_1$-norm at most $2$ has the
EJ property if and only if it has no odd-$K_4$ minor. Results of Conforti, Gerards and
Zambelli~\cite{ConGerZam07} and of Del Pia and Zambelli~\cite{dPZam09SIDMA} imply that every matrix
obtained from a totally unimodular matrix with two nonzero entries per row by multiplying some columns by
$2$ has the EJ property. Appa, Kotnyek, Papalamprou and Pitsoulis~\cite{AppKotPapPit07} showed that
integral binet matrices have the EJ property. Del Pia, Musitelli and Zambelli~\cite{dPMusZam18JCTB}
studied the EJ property for matrices obtained from $\{0,\pm1\}$-matrices with at most two nonzero entries
per column by multiplying some columns by $2$.

All matrices with the EJ property in these classes satisfy the following
condition~\cite{dPMusZam18JCTB}. An integral matrix $A$ is \emph{totally half-modular}
(\emph{THM}) if $S^{-1}$ is half-integral for every nonsingular square submatrix $S$ of $A$. Equivalently,
THM matrices are the integral $2$-regular matrices in the sense of Appa and Kotnyek~\cite{AppKot04}. In particular,
every entry of a THM matrix lies in $\{0,\pm1,\pm2\}$. For THM matrices the Chv\'atal closure has a simple
description by $\{0,\frac12\}$-combinations of the constraints~\cite{AppKot04,dP09thesis} (see
also~\cite{CapFis96} for $\{0,\frac12\}$-cuts).

A \emph{minor} of $A$ is any matrix obtained from $A$ by a sequence of the following
operations~\cite{dP09thesis,dPMusZam18JCTB}:
\begin{enumerate}[label=(\roman*)]
\item deleting, permuting, or multiplying by $-1$ rows or columns.
\item dividing by $2$ a row all of whose entries are even.
\item pivoting on an entry equal to $1$, that is, replacing
$\begin{pmatrix}1&\beta\\ \alpha&D\end{pmatrix}$ by $\begin{pmatrix}-1&\beta\\ \alpha&D-\alpha\beta\end{pmatrix}$, where $\alpha$ is a column
vector and $\beta$ is a row vector.
\end{enumerate}
A minor of $A$ is \emph{proper} if it has fewer rows plus columns than $A$. The EJ property and total
half-modularity are preserved under taking minors~\cite[Theorems~2.1 and~2.9]{dP09thesis} (see
also~\cite{GerSch86} for deletions and contractions). The two matrices
\[
A_3:=\begin{pmatrix}1&2&0\\1&2&2\\1&0&2\end{pmatrix},\qquad
A_4:=\begin{pmatrix}2&1&1&0\\2&1&0&1\\2&0&1&1\end{pmatrix}
\]
are THM~\cite[Section~1]{dPMusZam18JCTB} and do not have the EJ property~\cite[Observation~2.11]{dP09thesis}. Gerards and Schrijver
conjectured that they are the only obstructions. The conjecture was communicated by Gerards in 2007 and appears
in~\cite[Conjecture~2.1]{dP09thesis} and~\cite[Conjecture~1]{dPMusZam18JCTB}.

\begin{conjecture}[Gerards--Schrijver]\label{conj:GS}
A totally half-modular matrix has the EJ property if and only if it has no minor equal to $A_3$ or $A_4$.
\end{conjecture}

By the theorem of Gerards and Schrijver~\cite{GerSch86} (see~\cite[Section~2.2]{dP09thesis}), the
conjecture holds for the matrices whose rows have $\ell_1$-norm at most $2$. For these matrices, the relevant
obstruction is $A_4$, which corresponds to an odd-$K_4$ minor, although $A_4$ itself is not in this class. The
conjecture also holds for the THM matrices in the class of~\cite{dPMusZam18JCTB}, and there the relevant
obstruction is $A_3$.

\subsection*{The counterexample}

We show that Conjecture~\ref{conj:GS} is false. Let
\[
A_6:=\begin{pmatrix}
0&0&1&-1&0&0\\
0&1&1&0&2&0\\
0&1&0&1&0&2\\
1&0&1&0&2&0\\
1&0&0&1&0&2
\end{pmatrix}.
\]
As with $A_3$ and $A_4$, the subscript is the number of columns.

\begin{theorem}\label{thm:main}
The matrix $A_6$ is totally half-modular, does not have the EJ property, and has no minor equal to $A_3$ or $A_4$.
\end{theorem}

\begin{corollary}\label{cor:false}
Conjecture~\ref{conj:GS} is false.
\end{corollary}

A THM matrix is a \emph{minimal obstruction} if it does not have the EJ property, but every proper minor has it. The
matrices $A_3$ and $A_4$ are minimal obstructions~\cite[Section~1]{dPMusZam18JCTB}.

\begin{theorem}\label{thm:minimal}
The matrix $A_6$ is a minimal obstruction.
\end{theorem}

\subsection*{Proof ideas}

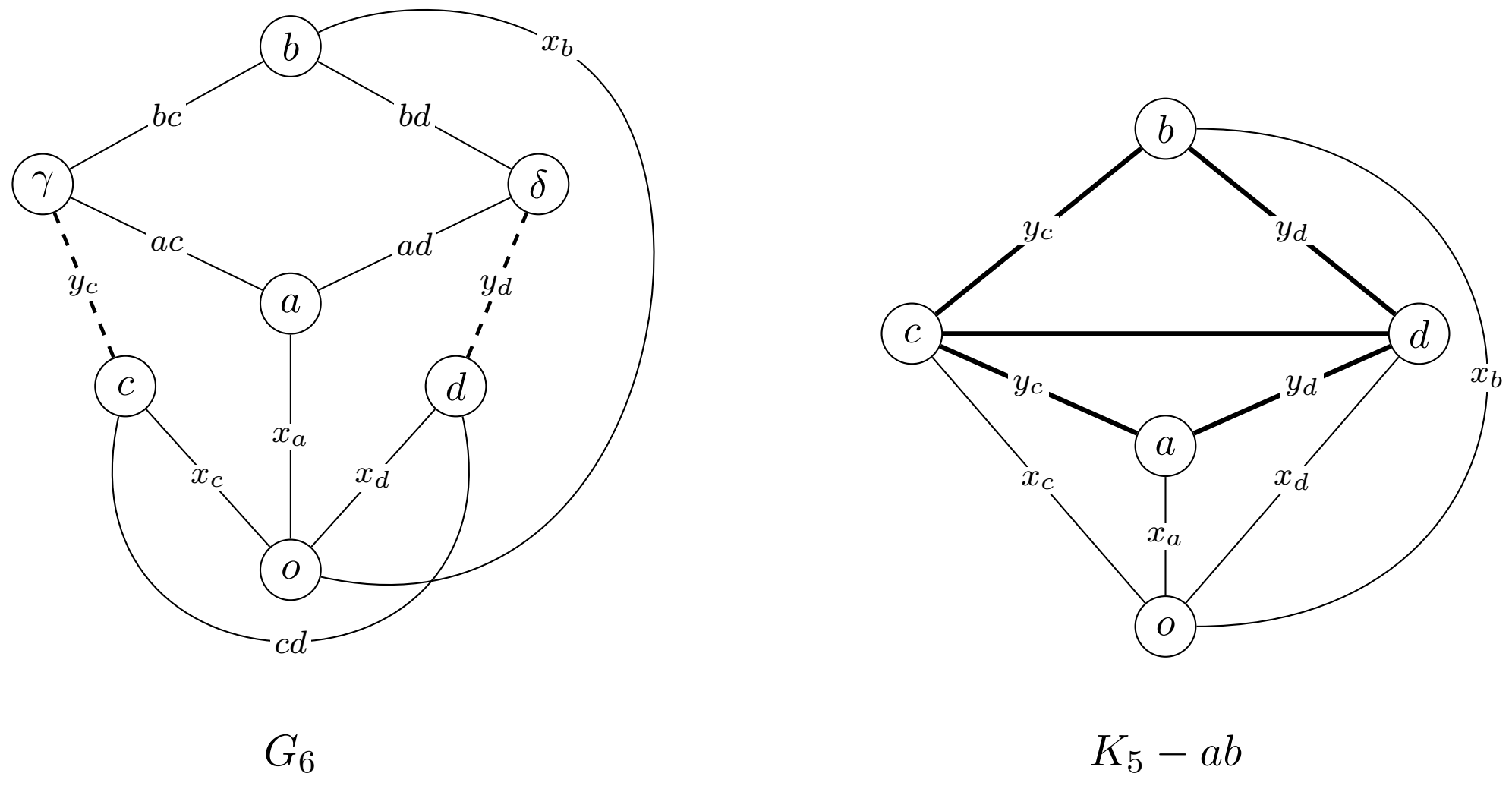
\begin{figure}[tb]
\centering
\begin{tikzpicture}[scale=0.8,
  v/.style={circle,draw,fill=white,inner sep=1pt,minimum size=15pt,font=\small},
  lab/.style={fill=white,inner sep=1pt,font=\scriptsize},
  lt/.style={dashed,line width=0.9pt}]
  \node[v] (o) at (0,-2.6) {$o$};
  \node[v] (c) at (-1.8,-0.6) {$c$};
  \node[v] (d) at (1.8,-0.6) {$d$};
  \node[v] (g) at (-2.7,1.6) {$\gamma$};
  \node[v] (h) at (2.7,1.6) {$\delta$};
  \node[v] (a) at (0,0.3) {$a$};
  \node[v] (b) at (0,3.1) {$b$};
  \draw (o) -- node[lab] {$x_a$} (a);
  \draw (o) -- node[lab] {$x_c$} (c);
  \draw (o) -- node[lab] {$x_d$} (d);
  \draw (o) .. controls (3.5,-3.4) and (4.6,0.6) .. (3.6,2.4) .. controls (2.9,3.6) and (1.2,3.7) .. node[lab,pos=0.25] {$x_b$} (b);
  \draw (a) -- node[lab] {$ac$} (g); \draw (a) -- node[lab] {$ad$} (h);
  \draw (b) -- node[lab] {$bc$} (g); \draw (b) -- node[lab] {$bd$} (h);
  \draw[lt] (c) -- node[lab] {$y_c$} (g);
  \draw[lt] (d) -- node[lab] {$y_d$} (h);
  \draw (c) .. controls (-2.6,-4.2) and (2.6,-4.2) .. node[lab,pos=0.5] {$cd$} (d);
  \node at (0,-4.6) {$G_6$};
\end{tikzpicture}
\hspace{1.2cm}
\begin{tikzpicture}[scale=0.85,
  v/.style={circle,draw,fill=white,inner sep=1pt,minimum size=15pt,font=\small},
  lab/.style={fill=white,inner sep=1pt,font=\scriptsize},
  row/.style={line width=1.2pt}]
  \node[v] (c) at (-2.6,0) {$c$};
  \node[v] (d) at (2.6,0) {$d$};
  \node[v] (o) at (0,-3.0) {$o$};
  \node[v] (a) at (0,-1.15) {$a$};
  \node[v] (b) at (0,2.1) {$b$};
  \draw (o) -- node[lab,pos=0.5] {$x_c$} (c);
  \draw (o) -- node[lab,pos=0.5] {$x_d$} (d);
  \draw (o) -- node[lab,pos=0.5] {$x_a$} (a);
  \draw (b) to[out=0,in=0,looseness=2.0] node[lab,pos=0.5] {$x_b$} (o);
  \draw[row] (c) -- (d);
  \draw[row] (a) -- node[lab,pos=0.55] {$y_c$} (c);
  \draw[row] (a) -- node[lab,pos=0.55] {$y_d$} (d);
  \draw[row] (b) -- node[lab,pos=0.5] {$y_c$} (c);
  \draw[row] (b) -- node[lab,pos=0.5] {$y_d$} (d);
  \node at (0,-4.3) {$K_5-ab$};
\end{tikzpicture}
\caption{The graphs $G_6$ and $K_5-ab$. In $G_6$, each edge is labeled by the row or the bound it encodes. Dashed edges
have half length. Contracting them gives $K_5-ab$, in which thick edges are rows and thin edges are bounds. Each row is
labeled by those of the variables $y_c$ and $y_d$ that have a nonzero coefficient in it.}
\label{fig:graphs}
\end{figure}

\paragraph{Total half-modularity (Section~\ref{sec:graphs}).} We denote the variables of the polyhedra
$P(A_6;b,c,l,u)$ by $x_a$, $x_b$, $x_c$, $x_d$, $y_c$, $y_d$, in the order of the columns of $A_6$, and we call $pq$ the row
of $A_6$ whose odd entries are the coefficients of $x_p$ and $x_q$. Every row has exactly two odd entries, so $A_6$ is
close to the transpose of the incidence matrix of a graph, with vertices $a,b,c,d$ for $x_a,x_b,x_c,x_d$ and an extra
vertex $o$ for the bounds. The columns of $y_c$ and $y_d$, whose entries are even, can be
absorbed by a change of variables. For example, $y_c$ appears only in the rows $ac$ and $bc$, with coefficient $2$, so in
terms of the new variable $x_c+2y_c$, which becomes a new vertex $\gamma$, these rows involve only two variables. Doing
the same for $y_d$, with a new vertex $\delta$, the rows and bounds of $A_6$ become the edges of the graph $G_6$ of Figure~\ref{fig:graphs} (Lemma~\ref{lem:graphs}). The bounds on $y_c$ and
$y_d$ become edges $c\gamma$ and $d\delta$ of half length. Total half-modularity follows because twice the vector of every edge is an integral combination of the vectors of
the edges of any spanning tree (Proposition~\ref{prop:THM}).

\paragraph{Failure of the EJ property (Section~\ref{sec:notEJ}).} We give integral data $b,c,l,u$ and a point $z$ that
lies in the Chv\'atal closure of $P:=P(A_6;b,c,l,u)$ but violates an inequality that is valid for $P_I$
(Proposition~\ref{prop:notEJ}). Since $A_6$ is THM, $z\in P'$ if $z$ satisfies every $\{0,\frac12\}$-cut, which adds
some of the inequalities defining $P$ whose normals sum to $0$ modulo $2$, divides by $2$, and rounds down. Which
inequalities combine depends only on their normals modulo $2$, and modulo $2$ the columns of $y_c$ and $y_d$ vanish. Contracting the two edges of
half length of $G_6$ gives $K_5$ minus the edge $ab$, which we denote by $K_5-ab$ (Figure~\ref{fig:graphs}). Its edges are
the rows of $A_6$, the row $pq$ being the edge $pq$, and the bounds on $x_a,x_b,x_c,x_d$, the bound on $x_p$ being the
edge $op$. At $z$, nine of the inequalities have slack $\frac13$ and all others have slack at least $1$. The nine
correspond to the edges of $K_5-ab$, and their normals are nonzero and pairwise distinct modulo $2$. So every
$\{0,\frac12\}$-cut adds inequalities with total slack at least $1$ at $z$, and it is satisfied at $z$.

\paragraph{Minors (Section~\ref{sec:nominor}).} The minor operations act on $G_6$ by deleting, contracting and halving
edges. So the rows and bounds of every minor of $A_6$ form a graph obtained from $G_6$ in this way. For $A_4$ they do
not form a graph, and for $A_3$ they form $K_4$ with two disjoint edges of half length, which cannot be obtained from
$G_6$ (Proposition~\ref{prop:nominor}).

\paragraph{Minimality (Section~\ref{sec:minimal}).} To show that every proper minor of $A_6$ has the EJ property, we reduce it to known classes. Changes of coordinates do not change the EJ property, and
they do not change the graph $G_6$ of Figure~\ref{fig:graphs}. So, up to the symmetries of the graph, every proper minor of $A_6$ is obtained, by operations that preserve the EJ property, from a matrix whose rows and bounds form the graph obtained from $G_6$ by deleting or contracting one edge, in one of nine cases (Lemma~\ref{lem:reduction} and Proposition~\ref{prop:minimal}). For five of
the nine graphs, the edges of half length have a common end, possibly after moving vertices of degree $2$ along their paths, and the
theorem of Del Pia and Zambelli applies (Lemma~\ref{lem:DZ}). For the other four,
we subdivide the edges of half length, and in one case split a vertex, and find a spanning tree $T$ such that the other edges, written in coordinates with respect to $T$, form a matrix
whose columns have $\ell_1$-norm at most $2$, so that the theorem of Edmonds and Johnson applies (Lemma~\ref{lem:trees}).

Theorem~\ref{thm:main} can also be checked by computer. The matrix $A_6$ is small, and it has finitely many
minors, since minors of THM matrices are THM, so their entries lie in $\{0,\pm1,\pm2\}$, and the operations (i)--(iii)
never increase the size. We prove Theorem~\ref{thm:main} by hand instead, because the proofs are short and show why $A_6$
is a counterexample, and because they develop the weighted graph $G_6$ and how the minor operations act on it, which is
also the basis of the proof of Theorem~\ref{thm:minimal}.

Section~\ref{sec:remarks} contains two further counterexamples and open questions.

\section{Total half-modularity}\label{sec:graphs}

In this section we write the rows and bounds of a matrix as vectors, and encode those of $A_6$ by the edges of the graph
$G_6$ of Figure~\ref{fig:graphs}. We then show that $A_6$ is THM.

As in Section~\ref{sec:intro}, the variables of $P(A_6;b,c,l,u)$ are $x_a,x_b,x_c,x_d,y_c,y_d$, and $pq$ is the row of
$A_6$ whose odd entries are the coefficients of $x_p$ and $x_q$.

\paragraph{Rows and bounds as vectors.} For an integral matrix $X$ with $n$ columns, let $E(X)$ be the family formed by the
rows of $X$ and the vectors $e_1^\top,\dots,e_n^\top$. Up to sign, these are the normals of the inequalities defining $P(X;b,c,l,u)$. Let
$M(X)$ be the matroid on $E(X)$ represented over $\Q$ by these vectors. Each operation (i)--(iii) of
Section~\ref{sec:intro} acts on $E(X)$ as follows, up to signs:
\begin{itemize}
\item Deleting a row deletes an element.
\item Deleting column $j$ \emph{contracts} $e_j^\top$, that is, it deletes $e_j^\top$ and the entry $j$ of every other
element.
\item Dividing a row by $2$ \emph{halves} an element of $E(X)$ that lies in $2\Z^n$.
\item Permutations and sign changes of rows reorder $E(X)$ and change signs. Permutations and sign changes of columns,
and the pivot on an entry $a_{ij}=1$, are unimodular changes of coordinates. The
pivot replaces the variable $x_j$ by $a_ix$, and it exchanges row $i$ and $e_j^\top$.
\end{itemize}
\paragraph{Weighted graphs.} Let $\Gamma$ be a graph whose vertex set contains a special vertex $o$, and in
which every edge $g$ has a weight $\kappa_g\in\{0,-1\}$. We call $g$ \emph{light} if $\kappa_g=-1$ and \emph{heavy}
otherwise. Let the vectors $f_p$, for the vertices $p\ne o$, form a basis of $\Q^n$, where $n$ is the number of these vertices,
and let $f_o$ be the zero vector. The family
$F(\Gamma)$ of $\Gamma$ is indexed by the edges of $\Gamma$, and the vector of an edge $g$ with ends $p,q$ is
\[
v_g:=2^{\kappa_g}(f_p-f_q),
\]
where the order of the ends does not matter, since we consider vectors only up to sign, as in the bullet points
above. A loop gives the zero vector. A family $F$ is \emph{encoded} by $\Gamma$ if the edges of $\Gamma$ are the elements of $F$
and, for some choice of the basis $f_p$, every edge $g$ has vector $\pm v_g$. Then the matroid of $F$ over $\Q$ is the cycle matroid of $\Gamma$, whose bases are the spanning trees of $\Gamma$, which
we regard as sets of edges.

For a family $F$ of vectors in $\Q^n$ and an element $g$ of $F$ with vector $v$, let $F\setminus g$ be $F$ without $g$, and let $F/g$ be the family
of the images of the vectors of $F\setminus g$ in $\Q^n/\Q v$. So deleting column $j$ of $X$ replaces $E(X)$ by
$E(X)/e_j^\top$. For a non-loop edge $g$ of $\Gamma$, let $\Gamma/g$ be obtained by contracting $g$, where the new vertex
is $o$ if $g$ has $o$ as an end.

\begin{lemma}\label{lem:ops}
Let $g$ be an edge of a weighted graph $\Gamma$.
\begin{enumerate}[label=(\alph*)]
\item $F(\Gamma)\setminus g=F(\Gamma-g)$.
\item If $g$ is not a loop, then $F(\Gamma)/g=F(\Gamma/g)$.
\item If $g$ is heavy, then replacing $v_g$ by $\frac12v_g$ gives $F(\Gamma')$, where $\Gamma'$ is $\Gamma$ with $g$
light.
\end{enumerate}
\end{lemma}

\begin{proof}
Part (a) is clear. For (b), let $p$ and $q$ be the ends of $g$. The
images of $f_p$ and $f_q$ in $\Q^n/\Q v_g$ are equal, and the images of the vectors $f_w$, for the vertices $w\ne o$ of
$\Gamma/g$, form a basis. If $q=o$, then the image of $f_p$ is $0$. Part (c) follows from $\frac12v_g=2^{-1}(f_p-f_q)$.
\end{proof}

Every row of $A_6$ has exactly two odd entries. The variable $y_c$ appears only in the rows $ac$ and $bc$, with
coefficient $2$, so in terms of $x_c+2y_c$ these rows involve only two variables, and similarly for $y_d$. This gives
the graph $G_6$ of Figure~\ref{fig:graphs}.

\begin{lemma}\label{lem:graphs}
The family $E(A_6)$ is encoded by the graph $G_6$ of Figure~\ref{fig:graphs}, where the rows $ac$, $bc$, $ad$, $bd$
and $cd$ are the edges $a\gamma$, $b\gamma$, $a\delta$, $b\delta$ and $cd$, the bound on $x_p$ is the edge $op$ for
$p\in\{a,b,c,d\}$, and the bounds on $y_c$ and $y_d$ are the light edges $c\gamma$ and $d\delta$.
\end{lemma}

\begin{proof}
Let $f_a:=-e_a^\top$, $f_b:=-e_b^\top$, $f_c:=e_c^\top$, $f_d:=e_d^\top$, $f_\gamma:=f_c+2e_{y_c}^\top$ and
$f_\delta:=f_d+2e_{y_d}^\top$, where $e_a,e_b,e_c,e_d,e_{y_c},e_{y_d}$ are the unit vectors of the variables
$x_a,x_b,x_c,x_d,y_c,y_d$. The rows $ac$, $bc$, $ad$ and
$bd$ are $f_\gamma-f_a$, $f_\gamma-f_b$, $f_\delta-f_a$ and $f_\delta-f_b$. The bounds on $x_a,\dots,x_d$ are
$\pm f_a,\dots,\pm f_d$, and those on $y_c$ and $y_d$ are $\frac12(f_\gamma-f_c)$ and $\frac12(f_\delta-f_d)$. The row
$cd$ is $f_c-f_d$.
\end{proof}

For a family $F$ of vectors, let $\Lambda(F)$ be the lattice generated by its vectors.

\begin{lemma}\label{lem:THMlattice}
An integral matrix $X$ is THM if $2v\in\Lambda(B)$ for every basis $B$ of $M(X)$ and every $v\in E(X)$.
\end{lemma}

\begin{proof}
Let $X_{R,C}$ be a nonsingular square submatrix of $X$, and let $J$ be the set of the other
columns. The rows of $X$ in $R$ and the vectors $e_j^\top$ with $j\in J$ form a basis $B$ of $M(X)$. Up to
permutations, the matrix with these rows is
\[
S=\begin{pmatrix}X_{R,C}&X_{R,J}\\0&I\end{pmatrix},\qquad\text{and}\qquad
S^{-1}=\begin{pmatrix}Y&-YX_{R,J}\\0&I\end{pmatrix}\quad\text{with}\quad Y=X_{R,C}^{-1}.
\]
By assumption, $2e_i^\top\in\Lambda(B)$ for every $i$, that is, $2e_i^\top S^{-1}$ is integral. So $S^{-1}$ is
half-integral, and so is $Y$.
\end{proof}

\begin{proposition}\label{prop:THM}
The matrix $A_6$ is THM.
\end{proposition}

\begin{proof}
Let $B$ be a basis of $M(A_6)$, that is, a spanning tree of $G_6$, and let $L$ be the lattice
generated by the vectors $f_p-f_q$ for the edges $pq$ of $B$. For an edge $g$ with ends $p,q$, we have
$f_p-f_q=2^{-\kappa_g}v_g$ with $\kappa_g\in\{0,-1\}$, so $L\subseteq\Lambda(B)$. For every vertex $r\ne o$, since
$f_o=0$, the vector $f_r$ is the sum of the vectors $\pm(f_p-f_q)$ over the edges $pq$ of the path of $B$ from $o$ to $r$, so $f_r\in L$. So $f_p-f_q\in L$ for every edge $g$ with ends $p,q$, and
$2v_g=2^{1+\kappa_g}(f_p-f_q)\in\Lambda(B)$. By Lemma~\ref{lem:THMlattice}, $A_6$ is THM.
\end{proof}

\section{Failure of the EJ property}\label{sec:notEJ}

In this section we describe the Chv\'atal closure of the polyhedra $P(A;b,c,l,u)$ for THM matrices $A$, and we show that
$A_6$ does not have the EJ property.

For an integral matrix $A\in\Z^{m\times n}$, we write the system defining $P(A;b,c,l,u)$ as $\tilde Ax\le h$, where $\tilde A$ is
the matrix
\begin{equation}\label{eq:M}
\begin{pmatrix}A\\-A\\I\\-I\end{pmatrix}.
\end{equation}
The rows of $\tilde A$ are called \emph{sides}. The sides $a_i$ and $-a_i$ from row $a_i$ of $A$
correspond to $a_ix\le c_i$ and $-a_ix\le-b_i$, and are the \emph{upper} and \emph{lower side} of row $i$.
The sides $e_j^\top$ and $-e_j^\top$ correspond to the bounds $x_j\le u_j$ and $x_j\ge l_j$. For $x\in\R^n$,
$s(x):=h-\tilde Ax$ is the vector of \emph{slacks}.

Let $k:=2m+2n$ be the number of sides. A vector $\mu\in\{0,1\}^k$ is
\emph{odd} if $\mu^\top \tilde A\equiv0\pmod 2$ and $\mu^\top h$ is odd. We identify $\mu$ with the set of sides
in its support.

The following description of the Chv\'atal closure is essentially due to Appa and Kotnyek~\cite{AppKot04}. For matrices whose
rows have $\ell_1$-norm at most $2$, it goes back to Gerards and Schrijver~\cite{GerSch86}.

\begin{lemma}[{Appa and Kotnyek~\cite{AppKot04}}]\label{lem:closure}
Let $A$ be a THM matrix, let $P=P(A;b,c,l,u)$ with integral data, and let $\tilde Ax\le h$ be the system defining $P$.
Then
\[
P'=\{x\in P:\ \mu^\top s(x)\ge1 \text{ for every odd } \mu\in\{0,1\}^k\}.
\]
\end{lemma}

\begin{proof}
We first show that $\tilde A$ is THM. A nonsingular square submatrix $S$ of $\tilde A$ does not contain a row together with its
negative. So, up to signs and permutations of rows and columns, $S=\begin{pmatrix}A_{R,C}&A_{R,J}\\0&I\end{pmatrix}$, where
$R$ is the set of the rows of $A$ that give a row of $S$, $J$ is the set of the columns $j$ such that $e_j^\top$ or
$-e_j^\top$ gives a row of $S$, and $C$ is the set of the other columns of $S$. Since $S$ is nonsingular, $A_{R,C}$ is a
nonsingular square submatrix of $A$. Hence $S^{-1}=\begin{pmatrix}Y&-YA_{R,J}\\0&I\end{pmatrix}$ with
$Y=A_{R,C}^{-1}$, which is half-integral.

We now prove the equality. By~\cite[Lemma~2.10]{dP09thesis} (see also~\cite{AppKot04}), applied to the THM matrix $\tilde A$,
the Chv\'atal closure $P'$ is the set of points of $P$ that satisfy $\lambda^\top \tilde Ax\le\lfloor\lambda^\top h\rfloor$ for
every $\lambda\in\{0,\frac12\}^k$ such that $\lambda^\top \tilde A$ is integral. For such a $\lambda$, let $\mu:=2\lambda$, so that
$\mu^\top \tilde A\equiv0\pmod 2$. If $\mu^\top h$ is even, the inequality is valid for $P$. Otherwise $\mu$ is odd, and the
inequality is $\frac12\mu^\top \tilde Ax\le\frac12(\mu^\top h-1)$, which is equivalent to $\mu^\top s(x)\ge1$. Every odd $\mu$
arises in this way, with $\lambda:=\frac12\mu$. Hence $P'$ is the right-hand side of the claimed equality.
\end{proof}

For an odd $\mu$, we call the inequality $\mu^\top s(x)\ge1$ a \emph{$\{0,\frac12\}$-cut}.

The next consequence of Lemma~\ref{lem:closure} is our tool to show that a point is in the Chv\'atal
closure.

\begin{lemma}\label{lem:third}
Let $A$ be THM, let $P=P(A;b,c,l,u)$ with integral data, and let $z\in P$. Assume that every side has slack
$\frac13$ or at least $1$ at $z$, and that the normals of the sides with slack $\frac13$ are nonzero and
pairwise distinct modulo $2$. Then $z\in P'$.
\end{lemma}

\begin{proof}
Let $\mu$ be odd. If the support of $\mu$ contains a side with slack at least $1$, then $\mu^\top s(z)\ge1$.
Otherwise, all sides in the support of $\mu$ have slack $\frac13$. Their normals sum to $0$ modulo $2$, and
they are nonzero and pairwise distinct modulo $2$. Hence the support of $\mu$ has at least three elements,
and $\mu^\top s(z)\ge3\cdot\frac13=1$. By Lemma~\ref{lem:closure}, $z\in P'$.
\end{proof}

Contracting the light edges of $G_6$ gives the graph $K_5-ab$ of Figure~\ref{fig:graphs}. In $K_5-ab$, the row $pq$ is the edge $pq$ and the bound on $x_p$ is the edge $op$. Since the
columns of $y_c$ and $y_d$ have only even entries, the edges of $K_5-ab$ describe the rows of $A_6$ and the bounds on $x_a,x_b,x_c,x_d$
modulo $2$. This suffices for $\{0,\frac12\}$-cuts, since the condition $\mu^\top\tilde A\equiv0\pmod 2$ in the
definition of an odd $\mu$ involves the normals of the sides only modulo $2$.

Let
\begin{equation}\label{eq:data}
\begin{aligned}
b&=(0,2,1,1,0), &\qquad c&=(2,4,3,3,2),\\
l&=(-1,0,-1,0,-1,-1), &\qquad u&=(1,2,1,2,2,2),\\
z&=(\tfrac23,\tfrac13,\tfrac23,\tfrac13,\tfrac23,\tfrac13).
\end{aligned}
\end{equation}

\begin{proposition}\label{prop:notEJ}
Let $b,c,l,u$ and $z$ be as in~\eqref{eq:data}. Then $z$ lies in the Chv\'atal closure of $P(A_6;b,c,l,u)$ but not in
its integer hull. In particular, $A_6$ does not have the EJ property.
\end{proposition}

\begin{proof}
We use the variables and the names of the rows introduced in Section~\ref{sec:graphs}, and let $P=P(A_6;b,c,l,u)$.

We first show that $z\in P'$. The five rows of $A_6$ take the values $\frac13,\frac73,\frac43,\frac83,\frac53$ at $z$.
Hence $z\in P$. Exactly nine sides have slack $\frac13$ at $z$, and all the other sides have slack at least $\frac43$.
The nine sides are the lower sides of the first three rows, the upper sides of the last two rows, and the bounds
$x_a\le1$, $x_b\ge0$, $x_c\le1$ and $x_d\ge0$. Since the coefficients of $y_c$ and $y_d$ are even, the normals of
these nine sides are, modulo $2$, the transposes of the vectors
\[
e_c+e_d,\ e_b+e_c,\ e_b+e_d,\ e_a+e_c,\ e_a+e_d,\ e_a,\ e_b,\ e_c,\ e_d,
\]
that is, the edges of $K_5-ab$. They are nonzero and pairwise distinct. By Proposition~\ref{prop:THM} and Lemma~\ref{lem:third}, $z\in P'$.

Next, we show that every $x\in P\cap\Z^6$ satisfies $x_a-x_b\le0$, and hence so does every $x\in P_I$. Since
$z_a-z_b=\frac13>0$, this gives $z\notin P_I$. Let $x\in P\cap\Z^6$ and assume that $x_a-x_b>0$. Since $x_a\le1$ and
$x_b\ge0$, we have $x_a=1$ and $x_b=0$. The lower sides of the rows $bc$ and $bd$ give
$x_c+2y_c\ge2$ and $x_d+2y_d\ge1$, and the upper sides of the rows $ac$ and $ad$ give $x_c+2y_c\le2$ and
$x_d+2y_d\le1$. So $x_c+2y_c=2$ and $x_d+2y_d=1$. Hence $x_c=2-2y_c$ is even and $x_d=1-2y_d$ is odd. Since $x_c\le1$
and $x_d\ge0$, we get $x_c\le0$ and $x_d\ge1$, which contradicts the lower side $x_c-x_d\ge0$ of the row $cd$.
\end{proof}

\section{No minor equal to \texorpdfstring{$A_3$}{A3} or \texorpdfstring{$A_4$}{A4}}\label{sec:nominor}

In this section we track how the operations (i)--(iii) change the graph $G_6$. For every matrix $X$ obtained from $A_6$
by a sequence of these operations, $E(X)$ is encoded by a graph obtained from $G_6$ by deleting edges, contracting
edges and lowering weights (Lemma~\ref{lem:lineage}). This excludes $A_4$, whose matroid is not graphic
(Lemma~\ref{lem:A4}), and the light edges of these graphs exclude $A_3$ (Lemma~\ref{lem:K4} and
Proposition~\ref{prop:nominor}).

Let $N(X)$ be the matroid on $E(X)$ represented over $\mathrm{GF}(2)$ by the reductions modulo $2$ of its elements.
An element of a matroid is a \emph{coloop} if it lies in every basis, or equivalently, if it is not in the span of the
other elements.

\begin{lemma}\label{lem:Ncoloop}
In a sequence of operations (i)--(iii), no element that is a coloop of $N(X)$ for the current matrix $X$ is deleted.
\end{lemma}

\begin{proof}
Only rows are deleted. The reductions modulo $2$ of the elements $e_1^\top,\dots,e_n^\top$ of $E(X)$ form a basis of
$\mathrm{GF}(2)^n$, so they form a basis of $N(X)$ that contains no row. Hence no row of $X$ is a coloop of $N(X)$.
\end{proof}

Let $E(X)$ be encoded by $\Gamma$. By Lemma~\ref{lem:ops}, deleting an element deletes it from $\Gamma$, halving an
element lowers its weight by $1$, and contracting an element that is not a loop contracts it in $\Gamma$.
Unimodular changes of coordinates do not change $\Gamma$. If a graph is obtained from $\Gamma$ by deleting edges,
contracting non-loop edges and lowering weights, each of its vertices $p$ is the image of a set $S(p)$ of vertices of
$\Gamma$. These sets are pairwise disjoint, each of them is connected by the contracted edges, and each edge of the new
graph is an edge of $\Gamma$ between two of them.

\begin{lemma}\label{lem:light}
Let $E(X)$ be encoded by $\Gamma$. The \emph{light graph} of $\Gamma$ has the same vertices as $\Gamma$ and the light edges as edges.
A \emph{light component} is a connected component of the light graph.
\begin{enumerate}[label=(\alph*)]
\item A light edge that is a bridge of the light graph is a coloop of $N(X)$.
\item A light edge that is not a loop does not lie in $2\Z^n$.
\item If a heavy edge that is not a loop lies in $2\Z^n$, then its ends are in the same light component.
\end{enumerate}
\end{lemma}

\begin{proof}
For a set $W$ of vertices not containing $o$, let $\psi_W$ be the linear function with $\psi_W(f_p)=1$ for $p\in W$
and $\psi_W(f_p)=0$ for the other vertices $p\ne o$. Then $\psi_W(f_p-f_q)\in\{0,\pm1\}$, and this value is $0$ if $p$
and $q$ are both in $W$ or both outside $W$. The elements of $E(X)$ generate $\Z^n$, so a linear function that is
integral on them is integral on $\Z^n$.

(a) Let $g$ be such an edge. Removing $g$ splits its light component into two parts, and we let $W$ be one that does
not contain $o$. The function $2\psi_W$ is odd on $v_g$ and even on every other element. So $v_g$ is not a sum of other
elements modulo $2$.

(b) Let $p\ne o$ be an end of the edge. The function $2\psi_{\{p\}}$ is integral on the elements and odd on the edge.

(c) Otherwise, let $W$ be the light component of one of the ends that does not contain $o$. The function $\psi_W$ is
integral on the elements, because every light edge has both ends in $W$ or both outside $W$, and it is odd on the edge.
\end{proof}

\begin{lemma}\label{lem:lineage}
Let $E(A)$ be encoded by $\Gamma$, and consider a sequence of operations (i)--(iii). Then, for every matrix $X$ in the
sequence, $E(X)$ is encoded by a graph obtained from $\Gamma$ by deleting edges, contracting non-loop edges and lowering
weights, no weight drops below $-1$, and the following hold.
\begin{enumerate}[label=(\alph*)]
\item For every light component $W$ of $\Gamma$, the images of the vertices of $W$ lie in one light component.
\item Every vertex of a light component with at least two vertices is the image of a vertex of a light component of
$\Gamma$ with at least two vertices.
\end{enumerate}
\end{lemma}

\begin{proof}
Dividing a zero row by $2$ does not change the matrix, so we may assume that only nonzero elements are halved, and they
are not loops. Deleting column $j$ contracts $e_j^\top\ne0$, which is not a loop. By Lemma~\ref{lem:ops}, the
first claim holds as long as no weight drops below $-1$. By Lemmas~\ref{lem:Ncoloop}
and~\ref{lem:light}(a), no light edge that is a bridge of the light graph is deleted,
and contractions only merge light components. By Lemma~\ref{lem:light}(b), light edges are not halved, so no weight
drops below $-1$, and by Lemma~\ref{lem:light}(c), a halved heavy edge joins two vertices of the same light component.
All claims follow by induction on the number of operations.
\end{proof}

\begin{figure}[tb]
\centering
\begin{tikzpicture}[scale=0.8,
  v/.style={circle,draw,fill=white,inner sep=1pt,minimum size=15pt,font=\small},
  lab/.style={fill=white,inner sep=1pt,font=\scriptsize},
  lt/.style={dashed,line width=0.9pt}]
  \node[v] (o) at (0,-1.6) {$o$};
  \node[v] (u) at (-2.4,1.4) {};
  \node[v] (t) at (2.4,1.4) {};
  \node[v] (s) at (0,0.4) {};
  \draw (o) -- node[lab] {$x_1$} (u);
  \draw[lt] (o) -- node[lab] {$x_2$} (s);
  \draw[lt] (u) -- node[lab] {$x_3$} (t);
  \draw (u) -- node[lab] {$\rho_1$} (s);
  \draw (s) -- node[lab] {$\rho_2$} (t);
  \draw (o) -- node[lab] {$\rho_3$} (t);
\end{tikzpicture}
\caption{A graph encoding $E(A_3)$, where $x_1,x_2,x_3$ are the variables and $\rho_1,\rho_2,\rho_3$ the rows of $A_3$.
Each edge is labeled by the row or the bound it encodes. Dashed edges are light.}
\label{fig:A3}
\end{figure}
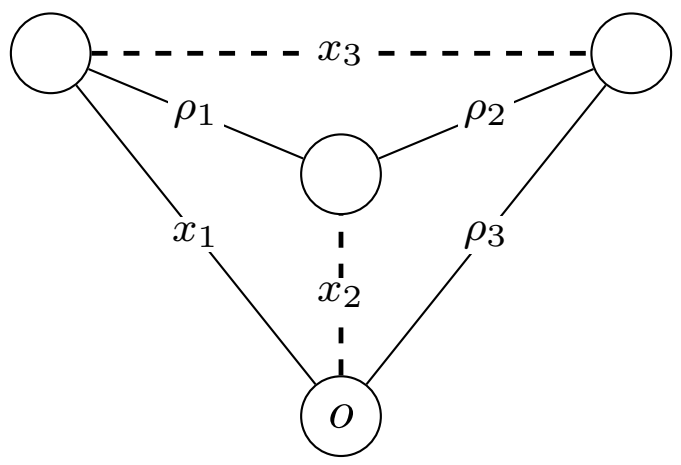

\begin{lemma}\label{lem:K4}
Let $E(A_3)$ be encoded by a graph $\Gamma$. Then $\Gamma$ is $K_4$, and its light edges form a perfect
matching.
\end{lemma}

\begin{proof}
The matroid $M(A_3)$ is the cycle matroid of $K_4$. Its $3$-element circuits are given by
\[
\rho_1=e_1^\top+2e_2^\top,\quad \rho_3=e_1^\top+2e_3^\top,\quad \rho_2=\rho_1+2e_3^\top,\quad \rho_2=\rho_3+2e_2^\top,
\]
where $\rho_1,\rho_2,\rho_3$ are the rows of $A_3$. Each of these relations is a linear dependency among three elements of $E(A_3)$ with
coefficients $\pm1$, $\pm1$ and $\pm2$, so exactly one coefficient is even.

The cycle matroid of $\Gamma$ is $M(A_3)$, which is connected, simple, of rank $3$ and with $6$ elements. So $\Gamma$ has
four vertices and is $K_4$. The three edges $g$ of a triangle satisfy $\sum\pm2^{-\kappa_g}v_g=0$. Up to a factor $2$,
this is the unique primitive integral relation of the triangle, and it has exactly one even coefficient if and only if
the triangle has exactly one light edge. Since $F(\Gamma)$ is $E(A_3)$ up to signs, these relations are, up to sign and a
factor $2$, the four relations above. So every triangle has
exactly one light edge, and the light edges form a perfect matching.
\end{proof}

\begin{lemma}\label{lem:A4}
The matroid $M(A_4)$ is not graphic.
\end{lemma}

\begin{proof}
The matroid $M(A_4)$ has $7$ elements and rank $4$. One checks directly that it is connected and has no circuit with at
most $3$ elements. (It is the dual of the non-Fano matroid.) A graph with cycle matroid $M(A_4)$ would be a simple graph
on $5$ vertices with $7$ edges and no triangle, but such a graph has at most $6$ edges.
\end{proof}

\begin{proposition}\label{prop:nominor}
The matrix $A_6$ has no minor equal to $A_3$ or $A_4$.
\end{proposition}

\begin{proof}
Suppose that $A'\in\{A_3,A_4\}$ is a minor of $A_6$, and fix a sequence of operations from $A_6$ to
$A'$. We first show that $A'=A_3$. By Lemmas~\ref{lem:graphs} and~\ref{lem:lineage}, $E(A')$ is encoded by a graph $\Gamma$ obtained from $G_6$ as in
Lemma~\ref{lem:lineage}. So $M(A')$ is graphic, and $A'=A_3$ by Lemma~\ref{lem:A4}.

Next, we identify the light components of $\Gamma$. By Lemma~\ref{lem:K4}, $\Gamma$ is $K_4$ and its light edges form a
perfect matching (see Figure~\ref{fig:A3}). So $\Gamma$ has two light components $W_1$ and $W_2$, each with two vertices, and the four heavy edges
join each vertex of $W_1$ to each vertex of $W_2$. These edges are heavy in $G_6$, since weights only decrease. Let $\bar
p$ denote the image of a vertex $p$ of $G_6$. The light components of $G_6$ with at least two vertices are
$\{c,\gamma\}$ and $\{d,\delta\}$. By Lemma~\ref{lem:lineage}(b), every vertex of $W_1\cup W_2$ is the image of one of
$c,\gamma,d,\delta$. By Lemma~\ref{lem:lineage}(a), $\bar c$ and $\bar\gamma$ lie in the same light component, and so do
$\bar d$ and $\bar\delta$. Since $W_1$ and $W_2$ both contain images, we may assume that $W_1=\{\bar c,\bar\gamma\}$ and
$W_2=\{\bar d,\bar\delta\}$.

Finally, we derive a contradiction. Let $S(p)$ be as above for the vertices $p$ of $\Gamma$. Since $\bar c\ne\bar\gamma$
and $\bar d\ne\bar\delta$, the
edges $c\gamma$ and $d\delta$ are not contracted, and the sets $S(p)$ are connected by heavy edges. The heavy edges at
$c$ go to $o$ and $d$, and those at $d$ go to $o$ and $c$. The graph $\Gamma$ has an edge between $\bar c$ and
$\bar\delta$. If $S(\bar c)=\{c\}$, this edge joins $c$ to $o$, so $o\in S(\bar\delta)$. Otherwise $S(\bar c)$
contains a heavy neighbor of $c$ other than $d$, so $o\in S(\bar c)$. Similarly, the edge between $\bar d$ and
$\bar\gamma$ gives $o\in S(\bar d)\cup S(\bar\gamma)$. These four sets are pairwise disjoint, a contradiction.
\end{proof}

\begin{proof}[Proof of Theorem~\ref{thm:main}]
The three claims are Propositions~\ref{prop:THM}, \ref{prop:notEJ} and~\ref{prop:nominor}.
\end{proof}

\begin{proof}[Proof of Corollary~\ref{cor:false}]
By Theorem~\ref{thm:main}, the matrix $A_6$ is THM, has no minor equal to $A_3$ or $A_4$, and does not
have the EJ property.
\end{proof}

\section{Minimality}\label{sec:minimal}

In this section we prove Theorem~\ref{thm:minimal}. In Section~\ref{sec:families} we show that the EJ property of a matrix $X$ depends
only on the family $E(X)$ up to changes of coordinates, in Section~\ref{sec:classes} we give two conditions on weighted
graphs under which their families have it, and in Section~\ref{sec:proofmin} we prove the theorem.

\subsection{Families and lattice bases}\label{sec:families}

Let $F$ be a finite family of vectors that spans a finite-dimensional vector space over $\Q$, such as $\Q^n$ or
$\Q^n/\Q v$, and recall that $\Lambda(F)$ is the
lattice generated by its vectors. A \emph{lattice basis} of $F$ is a subfamily of $F$ that is a basis of $\Lambda(F)$. If $B$ is a
lattice basis, writing the elements of $F$ outside $B$ in coordinates with respect to $B$ gives an integral matrix $X_B$,
with one row for each element outside $B$, and $E(X_B)$ is $F$ in these coordinates. For a matrix $X$ with $n$ columns, the
elements $e_1^\top,\dots,e_n^\top$ form a lattice basis $B$ of $E(X)$, and $X_B=X$. Two families $F$ and $F'$
\emph{correspond} under a linear isomorphism $\varphi$ if there is a bijection between them such that $\varphi$ maps the
vector of each element of $F$ to $\pm$ the vector of the corresponding element of $F'$.

The EJ property is not affected by repeated constraints.

\begin{lemma}\label{lem:trivial}
Let $A$ be an integral matrix, and let $\bar A$ be obtained from $A$ by adding a row equal to $\pm e_j^\top$ for some $j$,
or to $\pm$ a row of $A$. Then $\bar A$ has the EJ property if and only if $A$ has it.
\end{lemma}

\begin{proof}
Since $A$ is a minor of $\bar A$, the EJ property of $\bar A$ implies that of $A$~\cite[Theorem~2.1]{dP09thesis}.
Conversely, assume that $A$ has the EJ property, and let $\bar P=P(\bar A;b,c,l,u)$ with integral data. The two sides of the
added row can be merged with the bounds on $x_j$ or with the sides of the equal row. Hence $\bar P$ is either empty or
equal to a polyhedron $P(A;\cdot)$ with integral data. In the second case, $\bar P'=\bar P_I$ because $A$ has the EJ
property.
\end{proof}

\begin{lemma}\label{lem:families}
Let $F$ be a family as above that has a lattice basis.
\begin{enumerate}[label=(\alph*)]
\item Either all matrices $X_B$, for the lattice bases $B$ of $F$, have the EJ property, or none does. In the first case we say
that $F$ has the EJ property. Two families that correspond under a linear isomorphism both have the EJ property or
neither does.
\item Let $F$ have the EJ property, and let $B$ be a lattice basis. If $g\notin B$, then $F\setminus g$ has the EJ
property. If $g\in B$, then $F/g$ has the EJ property.
\item If two elements of $F$ are equal up to sign, then $F$ has the EJ property if and only if $F$ without one of them
has it.
\end{enumerate}
\end{lemma}

\begin{proof}
(a) Let $B$ and $B'$ be lattice bases. Since both are bases of $\Lambda(F)$, the change of coordinates from $B'$ to $B$ is
given by a unimodular matrix $U$, and $E(X_{B'})$ is $E(X_B)U$ up to signs and order. The system defining
$P(X_B;b,c,l,u)$ consists of two inequalities, with normals $v$ and $-v$, for each $v\in E(X_B)$, and the substitution
$x=Uy$ turns it into the system defining a polyhedron $P(X_{B'};b',c',l',u')$ with integral data. Every such polyhedron
arises in this way. The map $y\mapsto Uy$ preserves $\Z^n$, so it maps CG cuts to CG cuts, Chv\'atal closures to
Chv\'atal closures, and integer hulls to integer hulls. So $X_B$ has the EJ property if and only if $X_{B'}$ has it. A
linear isomorphism maps lattice bases to lattice bases and does not change the matrices $X_B$.

(b) If $g\notin B$, then $B$ is a lattice basis of $F\setminus g$, and $X_B$ for $F\setminus g$ is obtained from $X_B$ for
$F$ by deleting the row $g$. If $g\in B$, then the images of the elements of $B$ other than $g$ form a lattice basis of
$F/g$, and the corresponding matrix is obtained from $X_B$ by deleting the column of $g$. The EJ property is preserved
under taking minors~\cite[Theorem~2.1]{dP09thesis}.

(c) If one of the two elements is in a lattice basis $B$, the other is a row of $X_B$ equal to $\pm e_j^\top$. Otherwise,
both are rows of $X_B$. Apply Lemma~\ref{lem:trivial}.
\end{proof}

\subsection{Two classes of graphs with the EJ property}\label{sec:classes}

Let $\Gamma$ be a weighted graph as in Section~\ref{sec:graphs}. If $\Gamma$ is connected,
then $F(\Gamma)$ spans $\Q^n$. For a spanning tree $T$ of a connected graph and an edge $g\notin T$, the \emph{fundamental cycle} of $g$ is the
cycle formed by $g$ and the path of $T$ between the ends of $g$.

\begin{lemma}\label{lem:trees}
Let $\Gamma$ be a connected weighted graph without loops, and let $T$ be a spanning tree that contains all light edges.
\begin{enumerate}[label=(\alph*)]
\item $T$ is a lattice basis of $F(\Gamma)$.
\item If every heavy edge of $T$ lies on at most two fundamental cycles, and every light edge of $T$ lies on at most one,
then $F(\Gamma)$ has the EJ property.
\end{enumerate}
\end{lemma}

\begin{proof}
We first prove (a). Let $g\notin T$ be an edge with ends $p$ and $q$. It is heavy. For each edge $t$ of the path of $T$ between $p$ and $q$, the
vector $2^{-\kappa_t}v_t$ is $\pm(f_r-f_w)$, where $r$ and $w$ are the ends of $t$. Summing along the path,
\[
v_g=f_p-f_q=\sum_t\pm2^{-\kappa_t}v_t .
\]
So every element is an integral combination of the vectors of $T$, which are linearly independent. This proves (a).

We now prove (b). The column of $X_T$ indexed by $t\in T$ has the entry $\pm2^{-\kappa_t}$ in the rows of the edges whose fundamental cycle
contains $t$, and $0$ in the other rows. Under the assumption of (b), its $\ell_1$-norm is at most $2$. By the theorem of
Edmonds and Johnson~\cite{EdmJoh73}, $X_T$ has the EJ property.
\end{proof}

\begin{lemma}\label{lem:moves}
Let $\Gamma$ be a connected weighted graph, and let $p\ne o$ be a vertex with exactly two incident edges $g_1$ and $g_2$, with
other ends $q$ and $r$. Let $\Gamma'$ be obtained from $\Gamma$ by letting $g_1$ join $p$ to $r$ and $g_2$ join $p$ to $q$,
with the same weights. Then $F(\Gamma)$ and $F(\Gamma')$ correspond under a linear isomorphism in which each edge
corresponds to itself.
\end{lemma}

\begin{proof}
Let $f'_p:=f_q+f_r-f_p$, and $f'_w:=f_w$ for the other vertices $w\ne o$. These vectors form a basis of $\Q^n$. In
$\Gamma'$, the vector of $g_1$ is $2^{\kappa_{g_1}}(f'_p-f_r)=2^{\kappa_{g_1}}(f_q-f_p)$, and the vector of $g_2$ is
$2^{\kappa_{g_2}}(f'_p-f_q)=2^{\kappa_{g_2}}(f_r-f_p)$. These are the vectors of $g_1$ and $g_2$ in $\Gamma$ up to sign, and the
other vectors do not involve $f_p$.
\end{proof}

In Lemma~\ref{lem:moves}, the path $q\,p\,r$ is replaced by the path $r\,p\,q$, so the vertex $p$ moves along the
path and the two edges exchange their positions.

\begin{lemma}\label{lem:DZ}
Let $\Gamma$ be a connected weighted graph without loops whose light edges have a common end. Then $F(\Gamma)$ has the EJ
property.
\end{lemma}

\begin{proof}
We first reduce to the case where the common end is $o$, and we find a lattice basis. Let $r$ be a common end of the light edges, with $r:=o$ if there are no light edges.
Replacing each vector $f_p$ by $f_p-f_r$, including $f_o=0$, does not change the vectors
$f_p-f_q$, and $f_r$ becomes $0$. So $F(\Gamma)$ is also the family of $\Gamma$ with $r$ as the special vertex, and we may
assume that $r=o$. Let $\Gamma_0$ be obtained from $\Gamma$ by keeping one edge of each set of parallel edges with the
same weight. The light edges of $\Gamma_0$ form a star at $o$, so they lie in a spanning tree $T$ of $\Gamma_0$, which is a
lattice basis of $F(\Gamma_0)$ by Lemma~\ref{lem:trees}(a). Parallel edges with the same weight give the same vector up to
sign, so $T$ is also a lattice basis of $F(\Gamma)$, and by Lemma~\ref{lem:families}(c) we may replace $\Gamma$ by
$\Gamma_0$.

We now show that $F(\Gamma)$ has the EJ property. Let $\Gamma^+$ be obtained from $\Gamma$ by adding a heavy edge $op$ for every vertex $p\ne o$ not adjacent to $o$.
Its vector $f_p$ is an integral combination of the vectors of the edges of a path from $o$ to $p$, so
$\Lambda(F(\Gamma^+))=\Lambda(F(\Gamma))$, and $T$ is also a lattice basis of $F(\Gamma^+)$. By
Lemma~\ref{lem:families}(b), it suffices to show that $F(\Gamma^+)$ has the EJ property. For every vertex $p\ne o$,
choose an edge $g_p$ between $o$ and $p$, light if possible, and let $\sigma_p:=2$ if $g_p$ is light and $\sigma_p:=1$ otherwise.
The vector of $g_p$ is $f_p/\sigma_p$. The other edges between $o$ and $p$ are heavy, with vector $f_p=\sigma_p\cdot f_p/\sigma_p$, and
an edge between $p$ and $q$, with $p,q\ne o$, is heavy, with vector $\sigma_p\cdot f_p/\sigma_p-\sigma_q\cdot f_q/\sigma_q$. So the edges
$g_p$ form a lattice basis $B$, and $X_B=Y\Delta$, where $\Delta$ is the diagonal matrix with the entries $\sigma_p$, and every row
of $Y$ is $e_p^\top$ or $e_p^\top-e_q^\top$, with coordinates indexed by the vertices $p\ne o$. Adding to $Y$ a column with entry $-1$ in the rows with one nonzero
entry gives the transpose of the incidence matrix of a directed graph, which is totally unimodular and has two nonzero
entries per row. By the theorem of Del Pia and Zambelli~\cite[Corollary~2]{dPZam09SIDMA} (see~\cite{ConGerZam07} for the bipartite case), the matrix obtained from it by multiplying the columns by the
entries $\sigma_p$, and the new column by $1$, has the EJ property. Deleting the new column gives $X_B$.
\end{proof}

\subsection{Proof of Theorem~\ref{thm:minimal}}\label{sec:proofmin}

Let $F:=E(A_6)$, encoded by $G_6$ (Lemma~\ref{lem:graphs}).

\begin{lemma}\label{lem:lightbasis}
Every lattice basis of $F$ contains $c\gamma$ and $d\delta$.
\end{lemma}

\begin{proof}
The light edges $c\gamma$ and $d\delta$ are bridges of the light graph, so they are coloops of $N(A_6)$ by
Lemma~\ref{lem:light}(a). A lattice basis of $F$ is a basis of $\Z^6$, so its reduction modulo $2$ is a basis of
$N(A_6)$, and it contains every coloop.
\end{proof}

\begin{lemma}\label{lem:reduction}
Assume that $F(G_6-g)$ has the EJ property for every heavy edge $g$, and that $F(G_6/g)$ has the EJ property for every
edge $g$. Then $A_6$ is a minimal obstruction.
\end{lemma}

\begin{proof}
By Theorem~\ref{thm:main}, $A_6$ is THM and does not have the EJ property. Let $A'$ be a proper minor of $A_6$, and fix a
sequence of operations from $A_6$ to $A'$. It deletes a row or a column, and we let $A''$ be the matrix just before the
first deletion. Every row of $A_6$ has an odd entry, so no element of $E(A_6)$ lies in $2\Z^6$, and this does not change
under unimodular changes of coordinates. So no row is divided by $2$ before $A''$, and $E(A'')$ is $F$ after a unimodular
change of coordinates, which maps a lattice basis $B$ of $F$ to the elements $e_1^\top,\dots,e_6^\top$ of $E(A'')$. If the deletion deletes a row, it
deletes an element $g\notin B$, and $g$ is heavy by Lemma~\ref{lem:lightbasis}. The resulting matrix is $X_B$ for $F\setminus g=F(G_6-g)$ (Lemma~\ref{lem:ops}), after a change
of coordinates. Since $B$ is a lattice basis of $F\setminus g$, and $F(G_6-g)$ has the EJ property by assumption, this
matrix has the EJ property by Lemma~\ref{lem:families}(a). If the deletion deletes a column, it contracts an
element $g\in B$. The images of the other elements of $B$ form a lattice basis of $F/g=F(G_6/g)$, and the resulting
matrix is the matrix of $F/g$ for this lattice basis, after a change of coordinates. Since $F(G_6/g)$ has the EJ property
by assumption, this matrix has the EJ property by Lemma~\ref{lem:families}(a). The
EJ property is preserved by the remaining operations~\cite[Theorem~2.1]{dP09thesis}. So $A'$ has the EJ property.
\end{proof}

\begin{figure}[tb]
\centering
\begin{minipage}[b]{0.32\textwidth}\centering
\begin{tikzpicture}[scale=0.37,
  v/.style={circle,draw,fill=white,inner sep=0.5pt,minimum size=11pt,font=\scriptsize},
  lt/.style={dashed,line width=0.8pt},
  tr/.style={line width=1.4pt},
  ltr/.style={dashed,line width=1.4pt}]
  \node[v] (o) at (0,-2.6) {$o$}; \node[v] (a) at (0,0.3) {$a$}; \node[v] (b) at (0,3.1) {$b$};
  \node[v] (x) at (-2.25,0.5) {}; \node[v] (d) at (1.8,-0.6) {$d$}; \node[v] (h) at (2.7,1.6) {$\delta$};
  \draw (o) -- (a); \draw (o) -- (x); \draw (o) -- (d); \draw (o) .. controls (3.5,-3.4) and (4.6,0.6) .. (3.6,2.4) .. controls (2.9,3.6) and (1.2,3.7) .. (b);
  \draw (a) -- (x); \draw (a) -- (h); \draw (b) -- (x); \draw (b) -- (h);
  \draw[lt] (d) -- (h); \draw (x) .. controls (-2.8,-4.2) and (2.6,-4.2) .. (d);
\end{tikzpicture}\\[1pt]
{\small (1) $G_6/c\gamma$}
\end{minipage}\hfill
\begin{minipage}[b]{0.32\textwidth}\centering
\begin{tikzpicture}[scale=0.37,
  v/.style={circle,draw,fill=white,inner sep=0.5pt,minimum size=11pt,font=\scriptsize},
  lt/.style={dashed,line width=0.8pt},
  tr/.style={line width=1.4pt},
  ltr/.style={dashed,line width=1.4pt}]
  \node[v] (o) at (0,-2.6) {$o$}; \node[v] (a) at (0,0.5) {$a$}; \node[v] (b) at (0,3.1) {$b$};
  \node[v] (g) at (-2.7,1.6) {$\gamma$}; \node[v] (h) at (2.7,1.6) {$\delta$};
  \node[v] (y) at (-1.3,-0.9) {};
  \draw (o) -- (a); \draw (o) -- (y); \draw (o) to[bend left=30] (y); \draw (o) .. controls (3.5,-3.4) and (4.6,0.6) .. (3.6,2.4) .. controls (2.9,3.6) and (1.2,3.7) .. (b);
  \draw (a) -- (g); \draw (a) -- (h); \draw (b) -- (g); \draw (b) -- (h);
  \draw[lt] (y) -- (g); \draw[lt] (y) to[bend right=25] (h);
\end{tikzpicture}\\[1pt]
{\small (2) $G_6/cd$}
\end{minipage}\hfill
\begin{minipage}[b]{0.32\textwidth}\centering
\begin{tikzpicture}[scale=0.37,
  v/.style={circle,draw,fill=white,inner sep=0.5pt,minimum size=11pt,font=\scriptsize},
  lt/.style={dashed,line width=0.8pt},
  tr/.style={line width=1.4pt},
  ltr/.style={dashed,line width=1.4pt}]
  \node[v] (o) at (0,-2.6) {$o$}; \node[v] (a) at (0,0.3) {$a$}; \node[v] (b) at (0,3.1) {$b$};
  \node[v] (g) at (-2.7,1.6) {$\gamma$}; \node[v] (h) at (2.7,1.6) {$\delta$};
  \node[v] (d) at (1.8,-0.6) {$d$};
  \draw (o) -- (a); \draw (o) -- (d); \draw (o) to[bend right=35] (d); \draw (o) .. controls (3.5,-3.4) and (4.6,0.6) .. (3.6,2.4) .. controls (2.9,3.6) and (1.2,3.7) .. (b);
  \draw (a) -- (g); \draw (a) -- (h); \draw (b) -- (g); \draw (b) -- (h);
  \draw[lt] (o) -- (g); \draw[lt] (d) -- (h);
\end{tikzpicture}\\[1pt]
{\small (3) $G_6/oc$}
\end{minipage}\\[4pt]
\begin{minipage}[b]{0.32\textwidth}\centering
\begin{tikzpicture}[scale=0.37,
  v/.style={circle,draw,fill=white,inner sep=0.5pt,minimum size=11pt,font=\scriptsize},
  lt/.style={dashed,line width=0.8pt},
  tr/.style={line width=1.4pt},
  ltr/.style={dashed,line width=1.4pt}]
  \node[v] (o) at (0,-2.6) {$o$}; \node[v] (a) at (0,0.3) {$a$}; \node[v] (b) at (0,3.1) {$b$};
  \node[v] (g) at (-2.7,1.6) {$\gamma$}; \node[v] (h) at (2.7,1.6) {$\delta$};
  \node[v] (c) at (-1.8,-0.6) {$c$}; \node[v] (d) at (1.8,-0.6) {$d$};
  \draw (o) -- (a); \draw (o) -- (c); \draw (o) -- (d); \draw (o) .. controls (3.5,-3.4) and (4.6,0.6) .. (3.6,2.4) .. controls (2.9,3.6) and (1.2,3.7) .. (b);
  \draw (a) -- (g); \draw (a) -- (h); \draw (b) -- (g); \draw (b) -- (h);
  \draw[lt] (c) -- (g); \draw[lt] (d) -- (h);
\end{tikzpicture}\\[1pt]
{\small (4) $G_6-cd$}
\end{minipage}\hfill
\begin{minipage}[b]{0.32\textwidth}\centering
\begin{tikzpicture}[scale=0.37,
  v/.style={circle,draw,fill=white,inner sep=0.5pt,minimum size=11pt,font=\scriptsize},
  lt/.style={dashed,line width=0.8pt},
  tr/.style={line width=1.4pt},
  ltr/.style={dashed,line width=1.4pt}]
  \node[v] (o) at (0,-2.6) {$o$}; \node[v] (a) at (0,0.3) {$a$}; \node[v] (b) at (0,3.1) {$b$};
  \node[v] (g) at (-2.7,1.6) {$\gamma$}; \node[v] (h) at (2.7,1.6) {$\delta$};
  \node[v] (c) at (-1.8,-0.6) {$c$}; \node[v] (d) at (1.8,-0.6) {$d$};
  \draw (o) -- (a); \draw (o) -- (d); \draw (o) .. controls (3.5,-3.4) and (4.6,0.6) .. (3.6,2.4) .. controls (2.9,3.6) and (1.2,3.7) .. (b);
  \draw (a) -- (g); \draw (a) -- (h); \draw (b) -- (g); \draw (b) -- (h);
  \draw[lt] (c) -- (g); \draw[lt] (d) -- (h); \draw (c) .. controls (-2.6,-4.2) and (2.6,-4.2) .. (d);
\end{tikzpicture}\\[1pt]
{\small (5) $G_6-oc$}
\end{minipage}\hfill
\begin{minipage}[b]{0.32\textwidth}\centering
\begin{tikzpicture}[scale=0.37,
  v/.style={circle,draw,fill=white,inner sep=0.5pt,minimum size=11pt,font=\scriptsize},
  lt/.style={dashed,line width=0.8pt},
  tr/.style={line width=1.4pt},
  ltr/.style={dashed,line width=1.4pt}]
  \node[v] (o) at (0,-2.6) {$o$}; \node[v] (c) at (-1.8,-0.6) {$c$}; \node[v] (d) at (1.8,-0.6) {$d$};
  \node[v] (c1) at (-3.0,0.4) {$c'$}; \node[v] (d1) at (3.0,0.4) {$d'$}; \node[v] (b) at (0,3.3) {$b$};
  \node[v] (h) at (2.7,2.0) {$\delta$};
  \node[v] (g) at (-2.7,2.0) {$\gamma$}; \node[v] (a) at (0,0.5) {$a$};
  \draw[tr] (o) -- (c); \draw[tr] (o) -- (d); \draw[tr] (o) .. controls (4.0,-3.4) and (4.9,0.8) .. (3.9,2.7) .. controls (3.1,3.9) and (1.2,3.9) .. (b);
  \draw[tr] (a) -- (g); \draw (a) -- (h); \draw[tr] (b) -- (g); \draw[tr] (b) -- (h);
  \draw (c) -- (c1); \draw[ltr] (c1) -- (g); \draw (d) -- (d1); \draw[ltr] (d1) -- (h); \draw (c) .. controls (-2.6,-4.2) and (2.6,-4.2) .. (d);
\end{tikzpicture}\\[1pt]
{\small (6) $H-oa$}
\end{minipage}\\[4pt]
\begin{minipage}[b]{0.32\textwidth}\centering
\begin{tikzpicture}[scale=0.37,
  v/.style={circle,draw,fill=white,inner sep=0.5pt,minimum size=11pt,font=\scriptsize},
  lt/.style={dashed,line width=0.8pt},
  tr/.style={line width=1.4pt},
  ltr/.style={dashed,line width=1.4pt}]
  \node[v] (o) at (0,-2.6) {$o$}; \node[v] (c) at (-1.8,-0.6) {$c$}; \node[v] (d) at (1.8,-0.6) {$d$};
  \node[v] (c1) at (-3.0,0.4) {$c'$}; \node[v] (d1) at (3.0,0.4) {$d'$}; \node[v] (b) at (0,3.3) {$b$};
  \node[v] (h) at (2.7,2.0) {$\delta$};
  \node[v] (g) at (-2.7,2.0) {$\gamma$}; \node[v] (a) at (0,0.5) {$a$};
  \draw[tr] (o) -- (a); \draw[tr] (o) -- (c); \draw[tr] (o) -- (d); \draw[tr] (o) .. controls (4.0,-3.4) and (4.9,0.8) .. (3.9,2.7) .. controls (3.1,3.9) and (1.2,3.9) .. (b);
  \draw[tr] (a) -- (h); \draw[tr] (b) -- (g); \draw (b) -- (h);
  \draw (c) -- (c1); \draw[ltr] (c1) -- (g); \draw (d) -- (d1); \draw[ltr] (d1) -- (h); \draw (c) .. controls (-2.6,-4.2) and (2.6,-4.2) .. (d);
\end{tikzpicture}\\[1pt]
{\small (7) $H-a\gamma$}
\end{minipage}\hfill
\begin{minipage}[b]{0.32\textwidth}\centering
\begin{tikzpicture}[scale=0.37,
  v/.style={circle,draw,fill=white,inner sep=0.5pt,minimum size=11pt,font=\scriptsize},
  lt/.style={dashed,line width=0.8pt},
  tr/.style={line width=1.4pt},
  ltr/.style={dashed,line width=1.4pt}]
  \node[v] (o) at (0,-1.6) {$o$}; \node[v] (c) at (-1.8,-0.6) {$c$}; \node[v] (d) at (1.8,-0.6) {$d$};
  \node[v] (c1) at (-3.0,0.4) {$c'$}; \node[v] (d1) at (3.0,0.4) {$d'$}; \node[v] (b) at (0,3.3) {$b$};
  \node[v] (h) at (2.7,2.0) {$\delta$}; \node[v] (g) at (-2.7,2.0) {$\gamma$};
  \draw[tr] (o) -- (c); \draw[tr] (o) -- (d); \draw[tr] (o) -- (b);
  \draw[tr] (o) -- (g); \draw[tr] (o) -- (h); \draw (b) -- (g); \draw (b) -- (h);
  \draw (c) -- (c1); \draw[ltr] (c1) -- (g); \draw (d) -- (d1); \draw[ltr] (d1) -- (h); \draw (c) .. controls (-2.6,-4.2) and (2.6,-4.2) .. (d);
\end{tikzpicture}\\[1pt]
{\small (8) $H/oa$}
\end{minipage}\hfill
\begin{minipage}[b]{0.32\textwidth}\centering
\begin{tikzpicture}[scale=0.37,
  v/.style={circle,draw,fill=white,inner sep=0.5pt,minimum size=11pt,font=\scriptsize},
  lt/.style={dashed,line width=0.8pt},
  tr/.style={line width=1.4pt},
  ltr/.style={dashed,line width=1.4pt}]
  \node[v] (o) at (0,-2.6) {$o$}; \node[v] (c) at (-1.8,-0.6) {$c$}; \node[v] (d) at (1.8,-0.6) {$d$};
  \node[v] (c1) at (-3.0,0.4) {$c'$}; \node[v] (d1) at (3.0,0.4) {$d'$}; \node[v] (b) at (0,3.3) {$b$};
  \node[v] (h) at (2.7,2.0) {$\delta$};
  \node[v] (w1) at (-1.6,1.6) {$w_1$}; \node[v] (w2) at (0.5,1.6) {$w_2$};
  \draw[tr] (o) -- (w1); \draw[tr] (o) -- (c); \draw[tr] (o) -- (d); \draw[tr] (o) .. controls (4.0,-3.4) and (4.9,0.8) .. (3.9,2.7) .. controls (3.1,3.9) and (1.2,3.9) .. (b);
  \draw (w1) -- (w2); \draw[tr] (w2) -- (b); \draw (w2) -- (h); \draw[tr] (b) -- (h);
  \draw (c) -- (c1); \draw[ltr] (c1) -- (w1); \draw (d) -- (d1); \draw[ltr] (d1) -- (h); \draw (c) .. controls (-2.6,-4.2) and (2.6,-4.2) .. (d);
\end{tikzpicture}\\[1pt]
{\small (9) $H'$}
\end{minipage}
\caption{The cases of the proof of Proposition~\ref{prop:minimal}. Dashed edges are light. Panels (1)--(5) show the graphs of
the cases, where unlabeled vertices are obtained by contracting an edge. Panels (6)--(9) show the graphs $\hat\Gamma$, with the
spanning trees $T$ in bold.}
\label{fig:cases}
\end{figure}
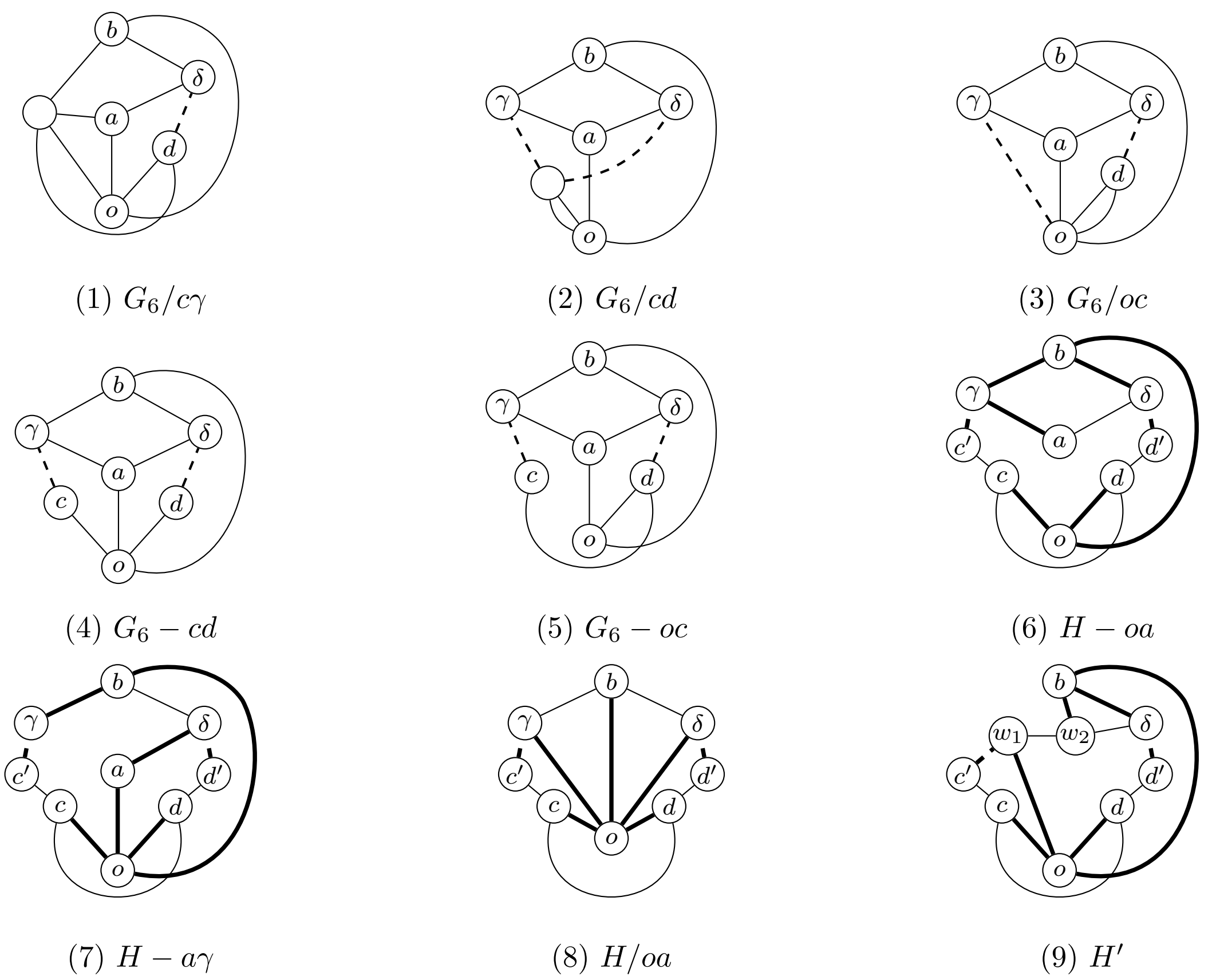

\begin{proposition}\label{prop:minimal}
For every heavy edge $g$ of $G_6$, the family $F(G_6-g)$ has the EJ property, and for every edge $g$ of $G_6$, so does
$F(G_6/g)$.
\end{proposition}

\begin{proof}
We first reduce to nine cases. The automorphisms of $G_6$ that exchange $a$ and $b$, or exchange $c$ with $d$ and $\gamma$ with $\delta$, fix $o$ and
the weights. They induce linear isomorphisms between the families, so by Lemma~\ref{lem:families}(a) it suffices to
consider one edge of each orbit, $\{oa,ob\}$, $\{oc,od\}$, $\{a\gamma,b\gamma,a\delta,b\delta\}$, $\{cd\}$ and
$\{c\gamma,d\delta\}$. Deleting the heavy edges $oa$, $oc$, $a\gamma$ and $cd$ and contracting the edges $oa$, $oc$,
$a\gamma$, $cd$ and $c\gamma$ gives nine cases: (1) $G_6/c\gamma$, (2) $G_6/cd$, (3) $G_6/oc$, (4) $G_6-cd$,
(5) $G_6-oc$, (6) $G_6-oa$, (7) $G_6-a\gamma$, (8) $G_6/oa$ and (9) $G_6/a\gamma$ (see Figure~\ref{fig:cases}). All these graphs are connected and have no loops.

Next, we treat cases (1)--(5) with Lemma~\ref{lem:DZ}.
\begin{enumerate}[label=(\arabic*)]
\item $G_6/c\gamma$ has one light edge, $d\delta$.
\item In $G_6/cd$, the light edges $c\gamma$ and $d\delta$ have the common end obtained from $c$ and $d$.
\item In $G_6/oc$, the vertex $c$ is identified with $o$, the light edge $c\gamma$ becomes $o\gamma$, and $cd$ and $od$
become parallel heavy edges with the same vector. We delete $cd$ by Lemma~\ref{lem:families}(c). Then $d$ has exactly two
incident edges, $od$ and $d\delta$, and Lemma~\ref{lem:moves} makes the light edge $d\delta$ join $d$ to $o$.
\item In $G_6-cd$, the vertex $c$ has exactly two incident edges, $oc$ and $c\gamma$, and so does $d$, with $od$ and
$d\delta$. Lemma~\ref{lem:moves} at $c$ and at $d$ makes both light edges end at $o$.
\item In $G_6-oc$, the vertex $c$ has exactly two incident edges, $cd$ and $c\gamma$. Lemma~\ref{lem:moves} makes the
light edge $c\gamma$ join $c$ to $d$, a common end with $d\delta$.
\end{enumerate}
In each case, Lemma~\ref{lem:DZ} applies to the final graph, and Lemmas~\ref{lem:moves} and~\ref{lem:families}(a), and
Lemma~\ref{lem:families}(c) in case (3), show that the family of the case has the EJ property.

Finally, we treat cases (6)--(9) with Lemma~\ref{lem:trees}. Let $H$ be obtained from $G_6$ by subdividing the light edges: $c\gamma$ is
replaced by a heavy edge $cc'$ and a light edge $c'\gamma$, and $d\delta$ by a heavy edge $dd'$ and a light edge
$d'\delta$, where $c'$ and $d'$ are new vertices. Contracting $cc'$ and $dd'$ in $H$ gives $G_6$. Let $H'$ be obtained
from $H/a\gamma$ by splitting the vertex $w$ obtained from $a$ and $\gamma$ into two vertices $w_1$, incident with the
images of $oa$ and $c'\gamma$, and $w_2$, incident with the images of $b\gamma$ and $a\delta$, joined by a new heavy edge
$w_1w_2$. Contracting $w_1w_2$ gives $H/a\gamma$. In each case, the graph of the
case is obtained by contractions from a graph $\hat\Gamma$, and we use a spanning tree $T$ of $\hat\Gamma$.
\begin{enumerate}[label=(\arabic*),start=6]
\item $G_6-oa=(H-oa)/\{cc',dd'\}$, and $T=\{ob,oc,od,a\gamma,b\gamma,b\delta,c'\gamma,d'\delta\}$.
\item $G_6-a\gamma=(H-a\gamma)/\{cc',dd'\}$, and $T=\{oa,ob,oc,od,b\gamma,a\delta,c'\gamma,d'\delta\}$.
\item $G_6/oa=(H/oa)/\{cc',dd'\}$, and $T=\{ob,oc,od,o\gamma,o\delta,c'\gamma,d'\delta\}$.
\item $G_6/a\gamma=H'/\{w_1w_2,cc',dd'\}$, and $T=\{ow_1,ob,oc,od,bw_2,b\delta,c'w_1,d'\delta\}$.
\end{enumerate}
In each case, $T$ contains the light edges of $\hat\Gamma$, and Figure~\ref{fig:cases} shows $\hat\Gamma$ with $T$ in bold. Each
light edge of $T$ lies only on the fundamental cycle of the edge $cc'$ or $dd'$ that subdivides it, since $c'$ and $d'$
have degree $2$ and $cc',dd'\notin T$. One checks from the figure that each heavy edge of $T$ lies on at most two
fundamental cycles. By Lemma~\ref{lem:trees}(b), $F(\hat\Gamma)$ has the EJ
property.

It remains to show that contracting $cc'$ and $dd'$, and also $w_1w_2$ in case (9), preserves the EJ property. The contracted edges and the light
edges form a forest: the paths $c\,c'\,\gamma$ and $d\,d'\,\delta$ in cases (6)--(8),
and the paths $c\,c'\,w_1\,w_2$ and $d\,d'\,\delta$ in case (9). So they lie in a spanning tree $T'$ of $\hat\Gamma$ that
contains all light edges, which is a lattice basis by Lemma~\ref{lem:trees}(a). We contract the edges one at a time. By
Lemma~\ref{lem:families}(b), each contraction preserves the EJ property, and the images of the remaining edges of $T'$
form a lattice basis. So the family of each case has the EJ property.
\end{proof}

\begin{proof}[Proof of Theorem~\ref{thm:minimal}]
Combine Lemma~\ref{lem:reduction} and Proposition~\ref{prop:minimal}.
\end{proof}

\section{Further counterexamples and open questions}\label{sec:remarks}

With similar arguments, we obtain two further counterexamples to Conjecture~\ref{conj:GS}. Let $A_6^+$ be obtained from $A_6$ by changing the entry $-1$ to $1$, and let $A_7$ be obtained from $A_6$ by changing the
last two entries of the first row to $2$ and $-2$ and adding the column $(2,0,0,0,0)^\top$:
\[
A_6^+:=\begin{pmatrix}
0&0&1&1&0&0\\
0&1&1&0&2&0\\
0&1&0&1&0&2\\
1&0&1&0&2&0\\
1&0&0&1&0&2
\end{pmatrix},
\qquad
A_7:=\begin{pmatrix}
0&0&1&-1&2&-2&2\\
0&1&1&0&2&0&0\\
0&1&0&1&0&2&0\\
1&0&1&0&2&0&0\\
1&0&0&1&0&2&0
\end{pmatrix}.
\]
Both matrices are THM, do not have the EJ property, and have no minor equal to $A_3$ or $A_4$. In the arguments, the
vector $f_c-f_d$ of the edge $cd$ of $G_6$ is replaced by $f_c+f_d$ for $A_6^+$, and the edge $cd$ of $G_6$ is replaced by
a path with one light edge for $A_7$. We do not include these proofs here.

Two matrices are \emph{equivalent} if they are related by pivots on entries $\pm1$, permutations and sign changes.
Conjecture~\ref{conj:GS} says that $A_3$ and $A_4$ are the only minimal obstructions, up to equivalence. Our
counterexamples show that this list is incomplete, and the natural question is how to complete it.

\begin{question}\label{q:list}
Does a THM matrix have the EJ property if and only if it has no minor equal to $A_3$, $A_4$, $A_6$, $A_6^+$ or $A_7$? More
generally, are there finitely many minimal obstructions, up to equivalence?
\end{question}

\subsection*{Acknowledgments}
The author thanks Giacomo Zambelli, his PhD advisor, for introducing him to the Gerards--Schrijver conjecture and
for the many hours they spent working on it together.

\bigskip
\begin{small}\sloppy
\paragraph{Funding.} A. Del Pia is partially supported by AFOSR Grant \mbox{FA9550-23-1-0433} and ONR Grant \mbox{N00014-25-1-2490}. Any opinions, findings, conclusions, or recommendations expressed in this material are those of the author and do not necessarily reflect the views of the Air Force Office of Scientific Research or the Office of Naval Research.

\paragraph{AI declaration.} AI tools (Claude, developed by Anthropic) were used for proofreading and suggesting improvements to the exposition. The author takes full responsibility for the content of the paper.
\end{small}

\bibliographystyle{amsplain}
\bibliography{refs}

\end{document}